\documentclass{article}
\usepackage{mathtools,amsmath,amssymb,amsthm,hyperref}

\usepackage[backend=bibtex]{biblatex}

\begin{document}
	\title{A question about points on an elliptic curve with prime denominator}
	\author{Simon L Rydin Myerson\thanks{Mathematics Insitute, University of Warwick; for contact details see \url{https://maths.fan}}}
	\date{7 July 2023\\Problem session of Journ\'ees Arithmetiques 2023\\Session chair: Michel Waldschmidt}
\maketitle
	\begin{abstract}
	Let \(E\) be an elliptic curve defined by a Weierstrass equation with integer coefficients. Any \(P\in E(\mathbb Q)\setminus\{O\}\) is of the form \(P=(\frac{x(P)}{z(P)^2},\frac{y(P)}{z(P)^3})\) where \(x(P),y(P)\in\mathbb Z\), \(z(P)\in \mathbb N\) and in addition both \(\gcd(x(P),z(P))=1\) and \(\gcd(y(P),z(P))=1\) hold. The question is: \textbf{how often is \(z(P)\) a prime number?}
	\end{abstract}
	
	\section{Preamble}
	
	This is an open problem submitted by the author and presented by Michel Waldschmidt to the JA2023 conference. I have done my best to include a number of different plausible strategies and questions for what I feel is a family of underexplored questions badly in want of new insights. Interested readers who may have questions or ideas to make progress are encouraged to  \href{https://maths.fan}{contact the author}; I would be  delighted to do what I can to help workers advance into this field.
	
	\section{What we expect}
	
	Let \(E(\mathbb Q) \simeq T\oplus \mathbb Z^R\) where \(T\) is the torsion subgroup. Let \(P_1,\dotsc,P_r\) generate the free part of  \(E(\mathbb Q) \) and write \(\vec{n}\cdot \vec{P}= n_1P_1+\dotsb+n_rP_r\).
	
	A general point of \(E(\mathbb Q) \) is of the form \(P=\vec{n}\cdot \vec{P}+S\) where \(\vec{n}\in \mathbb Z^r, S\in T\). For some constants \(0<c<C\) depending on \(E\), we have \(e^{c|\vec{n}|^2} < z(P) < e^{C|\vec{n}|^2} \). Heuristically the probability that \(z(P)\) is prime might be around \(1/\log z(P) \approx |\vec{n}|^{-2}\).
	
	Meanwhile the number of points \(P\) with \(N\leq |\vec{n}|\leq 2N\) is around \(N^{r}\).
	
	Let \(Z\) be a large parameter. Based on the above, we guess the total number of points with \(z(P)\) prime and \(z(P)\leq Z\) is around
	\begin{multline*}
	\#\{ P \in E(\mathbb Q) : z(P)\text{ prime}, z(P)\leq Z \}
	\\
	\approx
	\sum_{\substack{N=2^k, \, k\in \mathbb Z\\  e^{cN^2} <  Z } }N^{r-2}  \approx \begin{cases}\text{constant} &r\leq 1,
		\\
		\log \log Z &r=2,
		\\
		(\log Z)^{\frac{r}{2}-1} &r>2.
	\end{cases}
	\end{multline*}
	
	\section{What is known}
	\newtheorem{lem}{Lemma}
	\begin{lem}
		Suppose \(E\) is given by a minimal Weierstrass equation. If \(z(nP)\) is prime and \(n\) is sufficiently large, then \(n\) is prime.
		\end{lem}
		With a little work the condition that the Weierstrass equation is minimal can be removed (M. Verzobio, personal communication).
		\begin{proof}[Proof sketch.]
			Suppose that \(z(P)\) is prime and \(n=rs\) with \(r\geq s>1\). Then  because \( rP \equiv O \bmod z(rP) \) we must have \(nP\equiv O \bmod z(rP)\), which implies that \(z(rP)\) divides the denominator \(z(nP)\). If \(z(nP)\) is prime then \(z(rP)=z(nP)\) or \(z(rP)=1\). But \(r\geq \sqrt{n}\) and \(n\) is assumed large, so \(rP\) cannot be integral and \(z(rP)\neq 1\). Thus  \(z(rP)=z(nP)\).
			
			But if \(n\) is sufficiently large, then since the Weierstrass equation is minimal, \(z(nP)\) has a prime divisor which is coprime to \(z(rP)\), by a result of Silverman~\cite[pp232-233]{silverman1988}. This is impossible, so \(n\) must have been prime.
		\end{proof}
		If  \(E(\mathbb Q)\simeq \mathbb Z\), and \(Q\) is a generator of \(E(\mathbb Q)\), then
		\[
		\{ P \in E(\mathbb Q) : z(P)\leq Z \}
		\subseteq
		\{ nQ : n \ll \sqrt{\log Z}\}.
		\]
		From this and the lemma above we get:
		\newtheorem{cor}{Corollary}
		\begin{cor}
		 Suppose \(E\) is given by a minimal Weierstrass equation,
		 and that \(E(\mathbb Q)\simeq \mathbb Z\). Then
		 \[	\frac{\#\{ P \in E(\mathbb Q) : z(P)\text{ prime}, z(P)\leq Z \}}{
		 	\#\{ P \in E(\mathbb Q) : z(P)\leq Z \}
		 }
		 \ll
		 \frac{1}{\log \log Z}.
		 \]
		\end{cor}
		
		Actually there are some cases where much more is known. There are several results saying that for a suitable isogeny \(\phi:E'\to E\) and a suitable field \(K\), there are only finitely many points \(P\in E(\mathbb Q)\cap \phi(E'(K))\) such that \(z(P)\) is prime.

	As I understand it the proof idea is always something like this: Imagine for simplicity that \(K=\mathbb Q\). We try to show that the denominator \(z(\phi (Q))\) has to be divisible by more primes than the denominator \(z(Q)\), so if \(P=\phi(Q)\) has prime denominator then \(Q\) is integral and there are only finitely many such \(Q\). When \(K\neq \mathbb Q\) one has to look at the number of places at which \(Q\) has positive valuation, instead of primes dividing \(z(Q)\). I do not intend to provide a complete list of work of this kind, but I am thinking in particular of~\cite{einsiedlerEverestWard2001,everestKing2005,mahe2014}.
	\section{Specific questions}
	
	\begin{itemize}
		\item Show that Corollary 1 still holds if \(E\) has torsion, i.e.\@ if \(E(\mathbb Q)\simeq T\oplus \mathbb Z\).
		\item Do better than Corollary 1.
		That is take an elliptic curve for which 
				\[
				\#\{ P \in E(\mathbb Q) : z(P)\text{ prime}\}
				\]
		is \textbf{not} already known to be finite, and prove that
		\[
		\frac{\#\{ P \in E(\mathbb Q) : z(P)\text{ prime}, z(P)\leq Z \}}{
			\#\{ P \in E(\mathbb Q) : z(P)\leq Z \}/\log \log Z
		}\to 0 \qquad(Z\to \infty).
		\]
		\end{itemize}
		
		\section*{Acknowledgements}
		
		This question was inspired by my collaboration with Subham Bhakta, Dan Loughran, and Masahiro Nakahara.
		
		My subsequent exploration of these ideas has been carried out in the company of Subham Bhakta and Matteo Verzobio, to whom the author owes  many of the suggestions and hints in this document.
		
		I thank the organisers of Journ\'ees Arithm\'etiques for inviting me to participate, for their gracious hospitality and for their hard work to make the meeting run smoothly; I am grateful in particular to C\'ecile Dartyge, Youness Lamzouri and Paola Schneider.
		
		\printbibliography
		
\end{document}